\documentclass[12pt]{amsart}
\usepackage{amssymb, amsmath}

\newcommand{\bbc}{\mathbb{C}}
\newcommand{\bbp}{\mathbb{P}}
\newcommand{\bbq}{\mathbb{Q}}

\newcommand{\bbz}{\mathbb{Z}}
\newcommand{\Pic}{\mathrm{Pic}\,}
\newcommand{\prim}{\mathrm{prim}}

\newtheorem{thm}{Theorem}[section]

\newtheorem{prop}[thm]{Proposition}

\newtheorem{question}[thm]{Question}

\newtheorem{defnn}[thm]{Definition}

\newtheorem{remarkk}[thm]{Remark}
\newenvironment{remark}{\begin{remarkk} \em}{\end{remarkk}}
\newtheorem{examplee}[thm]{Example}

\title[Infinite multiplicity eigenvalues on elliptic curves]{Open conditions for infinite multiplicity eigenvalues on elliptic curves}
\author{Bo-Hae Im and Michael Larsen}
\date{\today}
\address{Department of Mathematics, University of Utah, Salt Lake City, Utah
 84112, USA}\email{im@math.utah.edu}
  \address{Department of Mathematics, Indiana University, Bloomington,
Indiana 47405, USA} \email{larsen@math.indiana.edu}

\begin{document}
\begin{abstract} Let $E$ be an elliptic curve defined over a
number field $K$.  We show that for each root of unity $\zeta$, the set
$\Sigma_\zeta$ of $\sigma\in\mathrm{Gal}(\overline{K}/K)$ such that $\zeta$ is an eigenvalue of 
infinite multiplicity for $\sigma$ acting on $E(\overline{K})\otimes\bbc$ has non-empty interior.

For the eigenvalue $-1$, we can show more: for any $\sigma$ in $\mathrm{Gal}(\overline{K}/K)$,
the multiplicity of the eigenvalue $-1$ is either $0$ or $\infty$.  It follows that  $\Sigma_{-1}$ 
is open.

\end{abstract}
\maketitle

\section{Introduction}
Let $K$ be a  number field, $\overline{K}$ an algebraic closure of
$K$, and $G_K := \mathrm{Gal}(\overline{K}/K)$ the absolute Galois group of
$\overline{K}$ over $K$. 
Let $E$ be an elliptic curve defined over $K$. 
There is a natural continuous action of $G_K$
on the countably infinite-dimensional complex vector space
$V_E := E(\overline{K})\otimes \bbc$.  The resulting representation decomposes as
a direct sum of finite-dimensional irreducible representations in each of which $G_K$
acts through a finite quotient group.

In particular, the action of every $\sigma\in G_K$ on $V_E$
is diagonalizable, with all eigenvalues roots of unity.   
In \cite{im2}, the first-named author showed that for \emph{generic} $\sigma$, every
root of unity appears as an eigenvalue of countably infinite multiplicity.  This is true both
in terms of measure and of Baire category.  However, there exist $\sigma$
for which the spectrum is quite different: 
trivially, the identity and complex conjugation elements; less trivially,
examples which can be constructed for an arbitrary set $S$ of primes, such that 
$\zeta$ is an eigenvalue if and only if every prime factor of its order lies in $S$.

Throughout this paper, we will write $\Sigma_\zeta$ for the subset of $G_K$ consisting
of elements $\sigma$ acting as $\zeta$ on an infinite-dimensional subspace of  
$V_E$ ($E$ and $K$ being fixed).
For $\zeta=1$, a good deal is known.  In \cite{im1}, it is proved that
whenever $1$ appears as an eigenvalue of $\sigma$
at all, we have $\sigma\in\Sigma_1$.  It follows that $\Sigma_1$
is open.  By \cite{im3}, when $K=\bbq$, $\Sigma_1$ is all of $G_K$, and 
quite possibly this may be true without restriction on $K$.
We have already observed that $\Sigma_\zeta\neq G_K$ for $\zeta\neq 1$.  We can still hope for
positive answers to the following progression of increasingly optimistic questions:

\begin{question}
\label{easy}
Does $\Sigma_\zeta$ have non-empty interior for all $\zeta$?
\end{question}

\begin{question}
\label{medium}
Is $\Sigma_\zeta$ open for all $\zeta$?
\end{question}

\begin{question}
\label{hard}
Do all eigenvalues of $\sigma$ acting on $V_E$ appear with infinite
multiplicity?
\end{question}

In this paper, we give an affirmative answer to Question~\ref{easy} for all $\zeta$
and an affirmative answer to all three questions for $\zeta = -1$.  

The difficulty in proving such theorems is that placing $\sigma$ in a basic open subset $U$ of
$G_K$ amounts to specifying the action of $\sigma$ on a finite Galois extension $L$ of $K$.
By the Mordell-Weil theorem, $E(L)\otimes\bbc$ cannot provide an infinite eigenspace for
$\zeta$.  Thus, the intersection of eigenspaces
$$\bigcap_{\sigma\in U} V_E^{\sigma-\zeta}$$
is finite-dimensional.  Thus, the behavior of a finite collection of rational points must be enough to guarantee the existence of infinitely many linearly independent points on the curve with specified $\sigma$-action.  

We would like to thank L.~Moret-Bailly and the referee for correcting versions of
Proposition~\ref{generic-nonzero} appearing in earlier drafts of this manuscript.

\section{Multiplicity of the eigenvalue -1}

In this section, we answer Questions~\ref{medium}~and~\ref{hard} for $\zeta=-1.$

\begin{prop} \label{prop:one} Let $E/K$ be an elliptic curve over $K$.
Suppose $-1$ is an eigenvalue of the action of $\sigma\in G_K$ on $V_E$.
Then the $-1$-eigenspace of $\sigma$ is infinite-dimensional.
\end{prop}

\begin{proof}   As $-1$ is an eigenvalue of $\sigma$ acting on $V_E$, it is an eigenvalue
of $\sigma$ acting on $E(\overline{K})\otimes\bbq$.  Clearing denominators, there exists
a non-torsion $P\in E(\overline{K})$ such that $\sigma(P) + P \in E(\overline{K})_\mathrm{tor}$.
Replacing $P$ by a suitable positive integral multiple, $\sigma(P) = -P$.

Let $y^2=f(x)$ be a fixed Weierstrass equation
of $E/K$.  Let $P=(\alpha,\sqrt{f(\alpha)})$.
As $\sigma(P)=-P$, we have
$\alpha\in \overline{K}^{\sigma}$ but
$\sigma(\sqrt{f(\alpha)})=-\sqrt{f(\alpha)}$ so $
\sqrt{f(\alpha)}\notin \overline{K}^{\sigma}$. Then,
$\sqrt{f(\alpha)}\notin K(\alpha)$, since $K(\alpha)\subseteq
\overline{K}^{\sigma}$

Note that $f(\alpha)\in K(\alpha)\subseteq \overline{K}^{\sigma}$.
Let $c=f(\alpha)\in K(\alpha)$. We still have $\sigma\in $
Gal$(\overline{K}/K(\alpha))$ and $\sigma(\sqrt{c})=-\sqrt{c}$.

Let $E'/K(\alpha)$ denote the twist $y^2=cf(x)$.
Then, $E'$ has a rational point
$P'=(\alpha, f(\alpha))$ over $K(\alpha)$. 
The $\overline{K}$-isomorphism $\phi\colon E\to E'$ mapping
$(x,y)\mapsto (x,\sqrt{f(\alpha)y})$ sends $P$ to $P'$, so $P'$ is of infinite order on $E'$.
By (\cite[Theorem 5.3]{im1}, $E'(\overline{K}^{\sigma})$
has infinite rank. Let $\{P_i'=(x_i,
\sqrt{cf(x_i)})\}_{i=1}^{\infty}$ be an infinite sequence of
linearly independent points of $E'$ generating the infinite
dimensional eigenspace of $1$ of $\sigma$ in
$E'(\overline{K})\otimes\bbc$. Then, $\sigma(x_i)=x_i$ and
$\sigma(\sqrt{f(x_i)})=-\sqrt{f(x_i)}$ for all $i$, since
$\sigma(\sqrt{c})=-\sqrt{c}$.

Let $P_i=\phi^{-1}(P_i')=(x_i,\sqrt{f(x_i)})$. These are points of
the given elliptic curve $E$ such that $\sigma(P_i)=-P_i$ for all
$i$, since $\sigma(x_i)=x_i$ and
$\sigma(\sqrt{f(x_i)})=-\sqrt{f(x_i)}$.

The points $P_i$ are linearly independent because the $P'_i$ are so.
Therefore,
$\{P_i\otimes 1\}_{i=1}^{\infty}$ generates an infinite
dimensional subspace of the $-1$-eigenspace of $\sigma$ on
$V_E$.  This completes the proof.
\end{proof}

\begin{thm}\label{thm:-1} Let $E/K$ be an elliptic curve over $K$. Then, $\Sigma_{-1}$ is open.
\end{thm}

\begin{proof} 
We have already seen that if 
$\sigma\in \Sigma_{-1}$, we can choose a point $P\in E(\overline{K})$ of infinite
order such that $\sigma(P)=-P$.  By Proposition~\ref{prop:one}, $\tau(P)=-P$
implies $\tau\in\Sigma_{-1}$.  It follows that $\Sigma_{-1}$ contains the open neighborhood
$\{\tau\in G_K\mid \tau(P)=\sigma(P)\}$ of $\sigma$.

\end{proof}

\begin{remark}
The same argument shows that Questions~\ref{medium}~and~\ref{hard} have an affirmative answer for $\zeta=\omega$ (resp. $\zeta=i$) when $E$ has complex multiplication
by $\bbz[\omega]$ (resp. $\bbz[i]$).
\end{remark}

\section{Interior Points}

In this section, we show that for every root of unity $\zeta$, the set $\Sigma_\zeta$ contains a
non-empty open subset.  We assume that the order of $\zeta$ is $n\ge 3$, the case $n=1$ having been treated in \cite{im1}, and the case $n=2$ in Theorem~\ref{thm:-1}.

Our strategy will be to find points $Q_i\in E(\overline{K})$ such that
the $\sigma$-orbit of $Q_i$ has length $n$.  For each such point $Q_i$, we set
\begin{equation}
\label{eigenvector}
R_i:=\sum\limits_{j=0}^{n-1}\sigma^{j}(Q_i)\otimes\zeta^{-j}
\end{equation}
and observe that $R_i$ is a $\zeta$-eigenvector of $\sigma$
provided that it is non-zero.  

We therefore begin with the following proposition:

\begin{prop}
\label{generic-nonzero}
Let $X$ be a Riemann surface of genus $g\ge 3$ with an automorphism $\sigma$ of order $n\ge 3$.
Then $X$ contains a non-empty open set $U$ such that $x\in U$ implies that
$$\sum_{i=0}^{n-1} [\sigma^i x]\otimes \zeta^{-i}\neq 0$$
in $\Pic X\otimes\bbc$.
\end{prop}

\begin{proof}
We can regard $X$ as the group of complex points of a non-singular projective curve whose
Picard scheme has complex locus $\Pic X$.  Then $\Pic X\otimes \bbz[\zeta]$ is the group of complex points of a group scheme whose identity component 
$\mathrm{Pic}^0X\otimes \bbz[\zeta]$ is isomorphic to the $\phi(n)$th power of 
the Jacobian variety of this curve. 
The action of $\sigma$ on $X$ defines an action on $\Pic X$,
and  the map $\psi\colon\Pic X\to \Pic X\otimes \bbz[\zeta]$ given by 
$$\psi(y) = \sum_{i=0}^{n-1} \sigma^i y\otimes \zeta^{-i}$$
then comes from a morphism of group schemes.  The image of $\psi$
actually lies in $\mathrm{Pic}^0X\otimes \bbz[\zeta]$, and its kernel $P_\zeta^0$
is Zariski-closed in $\Pic X$.

The set $P_\zeta$ of $y$ such that $\psi(y)$ maps to 0 in $\Pic X\otimes \bbc$
is the union of all translates of $P_\zeta^0$ by torsion points of $\Pic X$.
Applying Raynaud's theorem \cite{Ray} (i.e., the proof of the Manin-Mumford conjecture)
to the image of $X$ in
$\Pic X/P_\zeta^0$, the intersection $X\cap P_\zeta$
is finite whenever $\dim \Pic X/P_\zeta^0\ge 2$.  It therefore suffices to prove 
that the Lie algebra of $P_\zeta^0$ is a
subspace of the Lie algebra of $\Pic X$ of codimension $\ge 2$ or, equivalently, that
the rank of the map $\psi_*$ of Lie algebras is at least $2$.   We identify the Lie algebra of $\Pic X$
in the usual way \cite[Ch.~2,~\S6]{GH} with $H^1(X,{\mathcal O}_X) = H^{0,1}(X)$.
Likewise, the Lie algebra of $\Pic X\otimes \bbz[\zeta]$ is isomorphic to
$H^{0,1}(X)\otimes_{\bbz}\bbz[\zeta]$.  For every $k$ prime to $n$, there exists
a morphism 
$$\phi_k\colon H^{0,1}(X)\otimes_{\bbz}\bbz[\zeta]\to H^{0,1}(X)$$
obtained from the embedding of $\bbz[\zeta]$ into $\bbc$ mapping $\zeta$ to $\zeta^k$:
$$\phi_k(v\otimes \zeta^i) = \zeta^{ik}v.$$
The composition of this map with $\psi_*$
is $\sum_{i=0}^{n-1} \zeta^{-ik}\sigma^i$.

Let $H^{0,1}_{\prim}$ (resp. $H^1_{\prim}(X(\bbc),\bbc))$ denote the subspace of $H^{0,1}$
(resp. $H^1(X(\bbc),\bbc)$) spanned by eigenvectors of $\sigma$ whose eigenvalues are
primitive $n$th roots of unity.  
If $v$ is an eigenvector of $\sigma$ in $H^{0,1}$ whose eigenvalue is a primitive $n$th 
root of unity $\zeta^k$, then $\phi_k(\psi_*(v)) = nv\neq 0$, while $\phi_j(\phi_*(v)) = 0$ for
all $j\neq k$.  
It follows that $\ker\psi_*\cap H^{0,1}_{\prim} = \{0\}$, so the rank of $\psi_*$ is at least
$\dim H^{0,1}_{\prim}$.  The Hodge decomposition 
$$H^1(X(\bbc),\bbc) = H^{0,1} \oplus \overline{H^{0,1}}$$
implies 
$$\dim H^1_{\prim}(X(\bbc),\bbc) = 2\dim H^{0,1}_{\prim}.$$
It suffices, therefore, to prove $\dim H^1_{\prim}(X(\bbc),\bbc)\ge 4$.

Let $R_{\bbc}(G)$ denote the ring of complex (virtual) representations of $G$.
For any subgroup $H$ of $G:=\langle \sigma\rangle$, let $R_{G/H}$ denote the 
regular representation of $G/H$ regarded as an element of $R_{\bbc}(G)$, and let 
$I_H := R_G - R_{G/H}$.  In particular, $I_{\{1\}} = 0$.
Regarded as an element of $R_{\bbc}(G)$, the $G$-representation $H^1(X(\bbc),\bbc)$ is
$$2g+(2h-2)I_G+\sum_{[x]\in X/G}I_{\mathrm{Stab}_G(x)},$$
where $h$ is the genus of $X/G$, and $\mathrm{Stab}_G(x)$ is the stabilizer of any element of $X$
representing the $G$-orbit $[x]$.  
This is worked out in the case that $h=0$ in \cite[Prop.~2.2]{larsen}, but the method 
(in which the character of $H^1(X(\bbc),\bbc)$ as a representation of $G$ is deduced from
the Hurwitz formula and the Lefschetz trace formula) works in general.

The dimension of $H^1_{\prim}(X(\bbc),\bbc)$ is therefore 
$(2h-2+r)\phi(n)$, where $r$ is the number of ramification points of the
cover $X\to X/G$.  This is positive except in two cases: the cyclic cover $\bbp^1\to\bbp^1$
of degree $n$ (necessarily ramified over two points) and a degree $n$ isogeny  of elliptic curves; these
have genus $0$ and $1$ respectively.  Otherwise, it is at least $4$ unless 
$2h-2+r = 1$ and $\phi(n)=2$.  The triples $(h,r,n)$ for which this happens are
$(0,3,3)$, $(0,3,4)$, $(1,1,3)$, and $(1,1,4)$.  None of these is consistent with the condition
$g \ge 3$.

\end{proof}

\begin{thm}\label{thm:nth} Let $E/K$ be an elliptic curve over a number field $K$.
For each root of unity $\zeta$, there
exists a nonempty open subset $\Sigma_\zeta$ of $\mathrm{Gal}(\overline{K}/K)$
such that the multiplicity of the eigenvalue $\zeta$ for $\sigma\in \Sigma_\zeta$  acting on  $E(\overline{K})\otimes\bbc$ is infinite.
\end{thm}

\begin{proof} Let $\zeta$ be an $n$th root of unity.  
Let $\lambda_1,\lambda_2, \lambda_3,\infty$ be the
ramification points of a double cover $E\to \bbp^1$, and let $\lambda$ denote the
cross-ratio of $(\lambda_1,\lambda_2,\lambda_3,\infty)$.
Choose $a,b\in \overline{K}$ such that the ordered
quadruple ($a,b,\zeta a,\zeta b$) satisfies 
$$\dfrac{(\zeta a-a)(\zeta b-b)}{(\zeta b-a)(\zeta a-b)} = \lambda$$
This is always possible; for instance, setting $a=1$, we get a non-trivial quadratic equation
for $b$, and since $\lambda$ is not $1$ or $\infty$, we have $b,\zeta b\not\in\{a,\zeta a\}$.
Thus the elliptic curves
$$ X_i: y^2=(x-\zeta^{i-1}a)(x-\zeta^{i-1}b)(x-\zeta^ia)(x-\zeta^ib), \mbox{ for } i=1,\ldots, n.$$
all have the same $j$-invariant as $E$.

Let $L= K(a,b,\zeta)$.
Fix $q\in K$ such that $L(\sqrt[n]{q})$ is
a Galois $\bbz/n\bbz$-extension of $L$. 
We claim that $\Sigma_\zeta$ contains the open set
$$U_\zeta:=\{\sigma\in\mathrm{Gal}(\overline{K}/L)\mid \sigma(\sqrt[n]{q}) = \zeta\sqrt[n]{q}\}.$$

Let $M=L(\sqrt[n]{q})$.
For $N$ any number field containing $M$, let $C_N$ denote the affine curve over $N$
$$\mathrm{Spec}\;N[x,y_1,\ldots,y_n]/
(P_1(x,y_1),\ldots,P_n(x,y_n),y_1\cdots y_n - (x^n-a^n)(x^n-b^n))$$
where 
$$P_i(x,y) = y^2 - (x-\zeta^{i-1}a)(x-\zeta^i a)(x-\zeta^{i-1}b)(x-\zeta^i b).$$
Note that the equation $y_1\cdots y_n - (x^n-a^n)(x^n-b^n)=0$ merely selects one of the
two irreducible components of the 1-dimensional affine scheme cut out by the other equations.

Let $X$ denote the compact Riemann surface which is the compactification of
$C_N(\bbc)$.   By the Hurwitz genus formula, the genus of $X$ is $(n-2)2^{n-2}+1$,
which is $\ge 3$ since $n\ge 3$.
For any $n$-tuple $(k_1,\ldots,k_n)\in\{0,1\}^n$ with even sum, the map
\begin{equation}
\label{action}
(x,y_1,\ldots,y_n)\mapsto 
\bigl(\zeta x,(-1)^{k_1}\zeta^2 y_n,(-1)^{k_2}\zeta^2 y_1, (-1)^{k_3}\zeta^2 y_2,\ldots,
(-1)^{k_n}\zeta^2 y_{n-1}\bigr)
\end{equation}
defines an automorphism $\sigma$ of $C_N$ and therefore of $X$.  As the $k_i$ have even sum,
$\sigma$ is of order $n$.  If $x\in \sqrt[n]{q} L^*$
and $\sigma\in U_\zeta$, then
$\sigma(x) = \zeta x$, so
$$\sigma(y_i)^2 = \zeta^4 y_{i-1}^2,$$
and so there exists an $n$-tuple $(k_1,\ldots,k_n)$ with even coordinate 
sum such that $\sigma$ acts on 
$Q:=(x,y_1,\ldots,y_n)$ by (\ref{action}).  By Proposition~\ref{generic-nonzero}, for all but
finitely many values of $x$,
$$R:=\sum_{i=0}^{n-1} \sigma^i(Q)\otimes \zeta^{-i}$$
is a non-zero eigenvector of $\sigma$ with eigenvalue $\zeta$.

Assume now that $N$ is a finite Galois extension of $M$.
Consider the morphism from $C_N$ to the affine line over $M$ given by
$(x,y_1,\ldots,y_n)\mapsto x$.  This is a branched Galois cover with Galois group
$\mathrm{Gal}(N/M)\times (\bbz/2\bbz)^{n-1}$.  There exists a Hilbert set of values
$t\in M$ such that the geometric points lying over $x=\sqrt[n]{q}t$ in $C_M$ consists of a single 
$\mathrm{Gal}(\overline{K}/M)$-orbit or, equivalently, 
$\mathrm{Gal}(M(y_1,\ldots,y_n)/M)\cong (\bbz/2\bbz)^{n-1}$ and $M(y_1,\ldots,y_n)$ is linearly
disjoint from $N$ over $M$.
As a Hilbert set of a finite extension of $L$
always contains some Hilbert set of $L$ (\cite[Ch.~9,~Prop.~3.3]{l83}),
it follows that there exists $t\in L$ such that
setting $x=\sqrt[n]{q}t$, relative to $M$,
the extension $M(y_1,\ldots,y_n)$ is linearly disjoint from $N$
and has Galois group $(\bbz/2\bbz)^{n-1}$.

We can therefore iteratively construct a sequence $t_1,t_2,\ldots\in L^*$
such that the extensions
\begin{multline*}
\qquad
M_i:=M\biggl(\sqrt{(\sqrt[n]{q}t_i-a)(\sqrt[n]{q}t_i-b)(\sqrt[n]{q}t_i-\zeta a)(\sqrt[n]{q}t_i-\zeta b)},\ldots,
\hfill \\ \hfill
\sqrt{(\sqrt[n]{q}t_i-\zeta^{n-1}a)(\sqrt[n]{q}t_i-\zeta^{n-1}b)(\sqrt[n]{q}t_i-a)(\sqrt[n]{q}t_i-b)}\biggr)
\qquad \\
\end{multline*}
are all linearly disjoint over $M$.  Let $Q_i$ be a point with $x$-coordinate
$\sqrt[n]{q}t_i$, and 
$R_i$ the corresponding $\zeta$-eigenvector of $\sigma$ given by (\ref{eigenvector}).
We claim that the $R_i$ span a space of infinite dimension.  
The $Q_i$ do so by \cite[Lemma 3.12]{im1}, and as the
$\zeta^{-j}$ are linearly independent over $\bbq$, it follows that the $R_i$ do so as well.

\end{proof}

We conclude with a question that does not seem to be directly amenable to the methods of this
paper:
\begin{question}
Does the set $\bigcap_{\zeta\in \bbc^*_{\mathrm{tor}}} \Sigma_\zeta$ of elements
of $G_K$ having generic spectrum on $V_E$ always have an interior point?
\end{question}

\end{document}